\documentclass{article}

\usepackage[utf8]{inputenc}
\usepackage{amsfonts}       
\usepackage{nicefrac,fancyhdr,fullpage,hyperref}       
\usepackage{tikz}
\usepackage{subfigure}
\usepackage{mathrsfs}
\usepackage{amsmath, amssymb}
\usepackage{amsthm}
\newtheorem{theorem}{Theorem}[section]
\newtheorem{claim}{Claim}[theorem]

\newtheorem{proposition}{Proposition}[section]
\newtheorem{corollary}{Corollary}[theorem]

\newtheorem{lemma}{Lemma}[section]
\theoremstyle{definition}%
\newtheorem{definition}{Definition}[section]

\newcommand{\rna}[1]{\sigma^-{(#1)}}

\newcommand{\subgraph}[1]{[ #1 ]}
\newcommand{\floor}[2]{\left \lfloor \frac{#1}{#2} \right \rfloor}
\newcommand{\ceil}[2]{\left \lceil \frac{#1}{#2}\right \rceil}
\newcommand{\cfloor}[2]{\left\lceil\frac{1}{2}\floor{#1}{#2}\right\rceil}

\newcommand{\cceil}[2]{\left\lceil\frac{1}{2}\ceil{#1}{#2}\right\rceil}

\newcommand{\cfl}[2]{\left\lceil\frac{1}{2}(\floor{#1}{#2}-l)\right\rceil}
\newcommand{\ccl}[2]{\left\lceil\frac{1}{2}(\ceil{#1}{#2}+l)\right\rceil}
\newcommand{\offBra}[1]{B_{1#1}^{\Delta}}

\title{A study on certain bounds of the rna number and some characterizations of the parity signed graphs}
\author{Mohan Ramu\thanks{Department of Mathematics,CHRIST (Deemed to be University), Bengaluru, India. Email: mohan.r@res.christuniversity.in} \and 
Joseph Varghese Kureethara\thanks{Department of Mathematics,CHRIST (Deemed to be University), Bengaluru, India. Email: frjoseph@christuniversity.in}}

\date{}

\begin{document}
\maketitle
\abstract{For a given graph $G$, let $f:V(G)\to \{1,2,\ldots,n\}$ be a bijective mapping. For a given edge $uv \in E(G)$, $\sigma(uv)=+$, if $f(u)$ and $f(v)$ have the same parity and $\sigma(uv)=-$, if $f(u)$ and $f(v)$ have opposite parity. The resultant signed graph is called a parity signed graph and the mapping $\sigma$ is called a parity signature of $G$. Let us denote a parity signed graph $S=(G,\sigma)$ by $G_\sigma$. Let $E^-(G_\sigma)$ be a set of negative edges in a parity signed graph and let $Si(G)$ be the set of all parity signatures for the underlying graph $G$. We define the \textit{rna} number of $G$ as $\rna{G}=\min\{|E^-(G_\sigma)|:\sigma \in Si(G)\}$.	In this paper, we prove a non-trivial upper bound in the case of trees: $\rna{T}\leq \ceil{n}{2}$, where $T$ is a tree of order $n+1$. We have found families of trees whose \textit{rna} numbers are bounded above by $\ceil{\Delta}{2}$ and also we have shown that for any $i\leq \ceil{n}{2}$, there exists a tree $T$ (of order $n+1$) with $\rna{T}=i$. This paper gives a characterization of graphs with \textit{rna} number 1 in terms of its spanning trees and also a characterization of graphs with \textit{rna} number 2.}

\noindent \textbf{keywords}:Signed graph, parity signed graph, \textit{rna} number, tree
\section{Introduction}
\subsection{Parity Signed graphs}\label{sec1.1}
\indent All the graphs considered in this paper are simple and connected, and for the basic terminologies, we refer to \cite{harary1969graph, west1999introduction}. Let $G$ be a graph with $p$ vertices and $q$ edges. Let us denote the order of $G$ as $p(G)$ and $ G$'s size as $q(G)$. Let $E(G)$ be the set of all edges of $G$ and $V(G)$ be the set of all vertices of $G$. If $A$ and $B$ are two disjoint sets of vertices of $G$, then $E(A, B)$ is a set of edges between $A$ and $B$. We denote the order of a set $C$ by $|C|$.\\ 
\indent A \textit{signed} graph $S=(G,\sigma)$ is a simple graph in which $+$ or $-$ signs are assigned to each edge of $G$, where $\sigma: E(G)\to \{+,-\}$ and $\sigma$ is called a signature. Frank Harary introduced the concept of signed graphs (\cite{harary_notion_1953}), and since then, the domain of signed graphs has received great attention. This paper deals with a specific type of signed graph called \textit{parity signed graph} first defined by \cite{acharya_parity_2021}. For a given graph $G$, let $f:V(G)\to \{1,2,\ldots,n\}$ be a bijective mapping, and $\sigma(uv)=+$, if $f(u)$ and $f(v)$ have the same parity and $\sigma(uv)=-$, if $f(u)$ and $f(v)$ have opposite parity. The resultant signed graph is called a parity signed graph. Let us denote a parity signed graph $S=(G,\sigma)$ by $G_\sigma$. Let $E^-(G_\sigma)$ be a set of negative edges in a parity signed graph and let $Si(G)$ be the set of all parity signatures for the underlying graph $G$. We define the \textit{rna} number of $G$ as $\rna{G}=\min\{p(E^-(G_\sigma)):\sigma \in Si(G)\}$.\\
\indent \cite{acharya_characterizations_2021}  studied parity signed graphs further. \cite{sehrawat_rna_2021} studied the $rna$ number in the case of generalized Petersen graphs.  \cite{jin2021study} solved few open problems and proved a non-trivial upper bound for the \textit{rna} number. \\
Let us recall with the characterization theorem for parity signed graph.
\begin{theorem}\cite{acharya_characterizations_2021}
	A signed graph S is a parity signed graph if and only if its vertex set
	V(S) can be partitioned into two subsets $V_1(S)$ and $V_2(S)$ such that negative edges
	lie across $V_1(S)$ and $V_2(S)$ and $||V_1(S)|-|V_2(S)|| \leq 1$.
\end{theorem}
\subsection{Structure of the paper and notations}
In section \ref{sec2.1}, we study the \textit{rna} number for trees and obtain non-trivial upper bounds for the \textit{rna} number of the trees: For a tree $T$ of order $n$, we get $\rna{T}\leq \ceil{n-1}{2}$. Moreover, we study the classes of trees whose $rna$ numbers are either $\leq \ceil{\Delta}{2}$ or $>\ceil{\Delta}{2}$, where $\Delta$ is the maximum degree of a graph(in this case, tree). In section \ref{sec2.2}, we prove that a class of trees of order n can be partitioned into classes $R_i$, where $R_i$ consists of graphs with $rna$ number $i$, $1\leq i\leq \ceil{n-1}{2}$. We prove this by constructing a specific family of trees of order n, in which each member's $rna$ number takes a different value from the set $\{1,2,\ldots,\ceil{n-1}{2}\}$. Finally, in section \ref{sec2.3}, we characterize graphs with $rna$ number equal to 1 using the spanning trees of the graph. We also prove the characterization of graphs with $rna$ number equal to 2.\\
\noindent \textit{Notations.} \\
\indent Let $G$ be a simple, connected graph of order $n$ (order and size of a graph is already defined in section \ref{sec1.1}) and we shall denote $\Delta(G)=\Delta$ as the maximum degree of $G$.. Let $N(n)=\{1,2,\ldots,n\}$. Let $N_o(n)\subset N$ and $N_e(n)\subset N$ be the set of all odd integers and the set of all even integers, respectively, such that $N_o(n)\cup N_e(n)=N$. Suppose $f:V(G)\longrightarrow N(n)$ is a bijective function, then $\sigma_f$ is a function $\sigma_f:E(G)\rightarrow \{-1,1\}$ such that $\sigma_f(uv)=1$ if $f(u),f(v)$ are of the same parity and $\sigma_f(uv)=-1$ if $f(u),f(v)$ are of different parity.  \\
If $H\subseteq V(G)$, then the subgraph of $G$ induced by $H$ is denoted by $G[H]$ (also called as induced subgraph of $G$).\\
Let $S=(G,\sigma)$ be a parity signed graph and $A,B$ be the partition of $S$. Then we know that $||A|-|B||\leq 1$ and the negative edges in $S$ lies between $A$ and $B$. Let us call such partition $A,B$ of a parity signed 
graph as parity partition.

\section{Bounds for the \textit{rna} number of trees}\label{sec2.1}
A non-trivial upper bound was found for the \(rna\) number of any simple graph $G$, that is, $\rna{G}\leq \floor{2m+n}{4}$, where $m$ and $n$ are the number of edges and the number of vertices of $G$\cite[Theorem 2.7]{jin2021study}. We have found a non-trivial bound for the class of trees.
\begin{theorem}
	Let $T_{n+1}$ be a tree of order $n+1>1$. Then, $\sigma^{-}(T_{n+1}) \leq \left \lceil \frac{n}{2}\right \rceil$ 
\end{theorem}

\begin{proof}
	We use induction on the number of vertices of the tree $T_{n+1}$.\\
	For $1< n+1 \leq 3$, the statement is true. 
	Let us assume that the statement is true for trees with order less than or equal to $n>3$. Consider a tree $T_{n+1}$ and create $T^\prime$ by deleting two pendant vertices $u,v$ from $T_{n+1}$. Label $u$ and $v$ with odd and even numbers, respectively. By induction hypothesis, $\sigma^{-}(T^{\prime})\leq\left \lceil \frac{n-2}{2} \right \rceil$. Let $xu$ and $yv$ be the edges incident with $u$ and $v$, respectively, in $T_{n+1}$. If both $xu$ and $yv$ are negative edges in $T_{n+1}$, then $x$ is labelled by an even number and $y$ is labelled by an odd number, where $x\neq y$. By switching the labels of $u$ and $v$, $xu$ and $yv$ will become positive edges. When $x=y$, $u \text{ and }v$ are attached to one vertex $x$.\\Therefore, $u$ and $v$ can induce at most $1$ negative edge in $T_{n+1}$. This implies, $\sigma^{-}(T_{n+1})\leq \sigma^{-}(T^\prime)+1$.   
\end{proof}

Let $S_{1}, S_{2}$ be a partition of $V(T)$ such that $T[S_{1}]$ and $T[S_2]$ are connected subgraphs of $T$. Define \[ C(T)=\min\left\{||S_{1}|-|S_{2}||: S_{1},S_{2} \subset V(T) \right\}\] 
We can see that $C(K_{1,n})=n-1$ and $C(T_{n+1})\leq n-2$ for $T_{n+1}\neq K_{1,n}, \: n>2$.\\
Now let us define the following parameter: a generalization of $C(T)$. 
\[ C(T,a)=\min\left\{||S_{1}|-|S_{2}||: S_{1},S_{2} \subset V(T) \text{ and } |S_{1}|\leq a \right\}\] 

\begin{theorem}
	For a tree $T_{n+1}, \sigma^{-}(T_{n+1})\leq \left \lfloor \frac{C(T_{n+1})}{2}\right \rfloor +1$.
\end{theorem}

\begin{proof}
	$S_{1}\text{ and } S_{2}$ be a partition of $V(T_{n+1})$ such that $C(T_{n+1})=||S_1|-|S_2||$. Let us assume $|S_{1}|\leq|S_{2}| \text{ and hence, }|S_{1}|\leq \left \lfloor \frac{n+1}{2}\right\rfloor , |S_{2}|\geq \left\lceil \frac{n+1}{2}\right\rceil$. Let $uv$ be the bridge between  $T[S_{1}]$ and $T[S_2]$. Label $S_{1}$ by even numbers, label a connected subgraph (containing $v$) of $T[S_{2}]$ of order $\left\lceil \frac{n+1}{2}\right\rceil$ by odd numbers and the remaining vertices of $T[S_{2}]$ by even numbers, which implies that the number of negative edges is at most $\left \lfloor \frac{C(T_{n+1})}{2}\right \rfloor$+1. Thus $ \sigma^{-}(T_{n+1})\leq \left \lfloor \frac{C(T_{n+1})}{2}\right \rfloor +1$
\end{proof}

\begin{corollary}\label{rna1}
	$\sigma^{-}(T_{n+1})=1 \text{ if and only if } C(T_{n+1})\leq 1$
\end{corollary}

\begin{proof}
	Using Theorem 3.5 (\cite{acharya_characterizations_2021}), we infer that $\sigma^{-}(T_{n+1})=1$ if and only if $C(T_{n+1})\leq 1$.
\end{proof}

\begin{theorem}\label{OnlyStar}
	For $n>2,\:\sigma^{-}(T_{n+1})=\left \lceil \frac{n}{2} \right \rceil \text{ if and only if } T_{n+1}=K_{1,n}$
\end{theorem}

\begin{proof}
	
	We know that, $\sigma^{-}(K_{1,n})=\left \lceil \frac{n}{2} \right \rceil$. Hence, the sufficiency is proved. \\
	To prove the necessary condition, consider $T_{n+1}\neq K_{1,n}$ which implies, $0\leq C(T_{n+1})\leq n-2$. Let $S_{1},S_{2}$ be a partition of $V(T_{n+1})$ such that $C(T_{n+1})=||S_1|-|S_2||$ and assume $|S_{1}|<|S_{2}|$. Let $uv$ be the bridge connecting  $T[S_{1}]$ and $T[S_2]$. Label $S_{1}$ by even numbers, and label maximal connected subgraph(including the vertex $v$) of order $\left \lceil \frac{n+1}{2} \right \rceil$ in $T[S_2]$ by odd numbers and label the remaining vertices in $\subgraph{S_2}$ by even numbers. Let $|E^{-}(T[S_2])|$ be the number of negative edges in $T[S_2]$. Then, $|E^{-}(T[S_{2}])|\leq |S_{2}|-\left \lceil \frac{n+1}{2} \right \rceil$. If $|S_{2}|-\left \lceil \frac{n+1}{2} \right \rceil\ = \left \lceil \frac{n}{2} \right \rceil\ - 1$ then $|S_{2}|=n$
	which is a contradiction to the fact that $C(T_{n+1})\leq n-2$. Hence, $|S_{2}|-\left \lceil \frac{n+1}{2} \right \rceil\ \leq \left \lceil \frac{n}{2} \right \rceil\ - 2$. As $uv$ is a negative edge in $T_{n+1}, \sigma^{-}(T_{n+1})\leq|E^{-}(T[S_2])|+1 \leq  \left \lceil \frac{n}{2} \right \rceil - 1$.
	
\end{proof}

\begin{corollary}\label{rna 2}
	If $|V(T)|=n\leq 7$, then $\sigma^-(T)\leq 2 \text{ whenever }T\neq K_{1,{n-1}}$.
\end{corollary}

\begin{proposition}\label{propo}
	If $T$ is a tree of odd order, then $\sigma^-(T) \leq \sigma^-(T-v)$, where $v$ is a pendant vertex of $T$.
\end{proposition}

\begin{proof}
	Let $u$ be a vertex adjacent to $v$ in $T$. Let $V_1$ and $V_2$ be the parity partition of $V(T-v)$ such that there are $\sigma^-(T-v)$ edges between them and $\big| V_1 \big|=\big| V_2 \big |$. Assume that $u\in V_1$ which implies $\big| V_1 \cup \{v\} \big| > \big| V_2 \big |$ in $T$. Thus $\sigma^-(T) \leq \sigma^-(T-v)$.
\end{proof}
\subsection{Bounds for the rna number of certain classes of trees}\label{sec2.1.1}
In this section, $\Delta(T)=\Delta$ is the maximum degree of the tree $T$. From Theorem \ref{OnlyStar} only stars attain the non-trivial upper bound $\ceil{n}{2}$ and $\Delta(K_{1,n})=n$. Now, "Is $\rna{T}= \ceil{\Delta(T)}{2}?$" is a question we ask, and we answer this question. Trivial solutions are the classes of paths and stars. This question fails for many classes of graphs (as we will see), and so, we ask the following question, "Is $\rna{T}\leq \ceil{\Delta}{2}$?", which again is not true and we provide a class of graphs which violate this inequality. Therefore, the trees can be grouped into two classes, one whose \(rna\) number is at most $\ceil{\Delta}{2}$ and the rest of whose \(rna\) numbers are above $\ceil{\Delta}{2}$. In this section, we provide degree bounds for the \(rna\) number of the trees.
\begin{definition}[\cite{harary1969graph}]
	A \textit{branch} at a vertex $u$ of a tree $T$ is a maximal subtree of $T$ containing $u$ as a leaf vertex (or pendant vertex). 
\end{definition}
It is easy to see that there are $deg(u)$ branches at the vertex $u$ and an $i$th branch at a vertex $u$ is denoted by $Br(u,i)$, where $1\leq i\leq deg(u)$. Let $SBr(u,i)=Br(u,i)-u$, which is a subtree of $Br(u,i)$ obtained by deleting the vertex $u$ from $Br(u,i)$, be called as a \textit{secondary branch} at a vertex $u$ and $1\leq i\leq deg(u)$. 
We shall index the secondary branches at the vertex $u$ in the decreasing sequence of their order, i.e., if $a>b$, then $p(\offBra{b})\geq p(\offBra{a})$. Therefore, $p(\offBra{i})\geq p(\offBra{i+1})$, $1\leq i\leq \deg(u)-1$.

\begin{lemma}\label{lemma2.1}
	Suppose $n$ is a non-zero positive integer. Then the following properties are true.
	\begin{enumerate}
		\item $2(\cceil{n}{2}-1)=\begin{cases}
			\frac{n}{2}-2 & \text{if }n\equiv 0 \pmod 4 \\
			\floor{n}{2} & \text{if }n \equiv 1 \pmod 4\\
			\floor{n}{2}-1 & \text{if }n \equiv r \pmod 4, r=2\text{ or }3
		\end{cases}$
		\item $\floor{n}{2}-1\geq \cceil{n}{2}$ if $n\geq 6$.  
	\end{enumerate}
\end{lemma}
\begin{proof}
	\textbf{(1)}: If $n\equiv 0 \pmod 4$, then $2(\cceil{n}{2}-1)=2(\frac{n}{4}-1)=\frac{n}{2}-2$.\\
	Let $n \equiv 1 \pmod 4$. Then $n$ is odd. If $4|(n+1)$, then $n \equiv 3 \pmod 4$. Therefore, $4\nmid (n+1)$ and  $2(\cceil{n}{2}-1)=2(\left \lceil\frac{1}{2}(\frac{n+1}{2}) \right\rceil-1)=2(\frac{1}{2}(\frac{n+1}{2}+1)-1)=\floor{n}{2}$.\\
	Let $n \equiv 2 \pmod 4$. Then $n$ is even and $4\nmid n$. This implies $2(\cceil{n}{2}-1)=2(\left \lceil\frac{1}{2}\frac{n}{2} \right\rceil-1)=2(\frac{1}{2}(\frac{n}{2}+1)-1)=\frac{n}{2}-1$.\\
	Let $n \equiv 3 \pmod 4$. Then $n$ is odd, and $4|(n+1)$. This implies $2(\cceil{n}{2}-1)=2(\left \lceil\frac{1}{2}(\frac{n+1}{2}) \right\rceil-1)=2(\frac{1}{2}(\frac{n+1}{2})-1)=\floor{n}{2}-1$.
	
	\medskip
	\noindent\textbf{(2)}: Given $n\geq 6$. We know that $$\left\lceil\frac{1}{2}\left\lceil\frac{n}{2}\right\rceil\right\rceil=\begin{cases}
		\frac{1}{2}(\frac{n}{2}-2)+1 & \text{if } n \equiv 0 \pmod 4\\
		\frac{1}{2}\floor{n}{2}+1 & \text{if } n \equiv 1 \pmod 4\\
		\frac{1}{2}(\floor{n}{2}-1)+1 & \text{if } n\equiv r \pmod 4, r=2 \text{ or }3 
	\end{cases}$$
	If $n\equiv k \pmod 4$, then $\floor{n}{2}-1\geq \cceil{n}{2}$, for each $0\leq k\leq 3$.
\end{proof}

\begin{theorem}\label{upperboundThm}
	Let $T$ be a tree of order $n$ and $u$ be a vertex (of $T$) of maximum degree $\Delta$. Let $\Delta=\ceil{n}{2}+l,l\geq 0$ and we shall assume that $p(\offBra{j})\geq p(\offBra{j+1}), 1\leq j<\Delta-1$. Then $\rna{T}\leq \ceil{\Delta}{2}$ under each of the following conditions.
	
	\begin{enumerate}
		
		\item $p(\offBra{1})=2$ and $l=0$, for $n\geq 6$.
		\item $\cfl{n}{2}\leq p(\offBra{1})\leq \floor{n}{2}-l$, where $0\leq l\leq \floor{n}{2}-1$.
	\end{enumerate}
\end{theorem}

\begin{proof}
	
	It is given that $\Delta(T)=\ceil{n}{2}+l,l\geq 0$ and let $d=\ceil{\Delta}{2}$. 
	
	\medskip
	\noindent \textbf{(1)}:
	For $n\geq 6$, we know that $\floor{n}{2}-1\geq \cceil{n}{2}$ (lemma \ref{lemma2.1}) and $p(SBr(u,1))=2$, which implies $p(\offBra{j})=2$ for all $j\leq \cceil{n}{2}$.
	There are three subcases.\\
	\textbf{Case 1.1}: If $4|n$, then $\sum_{j=1}^{d}p(\offBra{j})=\frac{n}{2}$. By labelling the vertices of $\offBra{j},1\leq j\leq d$ with even vertices and the rest of the vertices of $T$ with odd vertices, we get $\rna{T}\leq \frac{\Delta}{2}$. \\
	\textbf{Case 1.2}: If $n\equiv r \pmod 4, r=2 \text{ or }3$, then $\sum_{j=1}^{d-1}p(\offBra{j})=\floor{n}{2}-1$. By labelling the vertices of $\offBra{j}, 1\leq j \leq d-1,$ and a pendant vertex of $T$ in $\offBra{d}$ with even integers, and rest of the vertices of $T$ with odd integers, we get $\rna{T}\leq \ceil{\Delta}{2}$.  \\
	\textbf{Case 1.3}: If $n \equiv 1 \pmod 4$, then $\sum_{j=1}^{d-1}p(\offBra{j})=\floor{n}{2}$. Therefore, $\rna{T}\leq \ceil{\Delta}{2}$.
	
	\medskip
	\noindent\textbf{(2)}: There are two cases.\\
	\textbf{Case 2.1}: 
	Let $l\geq 0$ and $$\left\lfloor\frac{n}{2}\right\rfloor-l\geq p(\offBra{1})\geq\begin{cases} 
		\left\lceil{\frac{1}{2}(\frac{n}{2}-l)}\right\rceil+1 &, \text{ if }2|(\frac{n}{2}-l)\text{ \& }2|n\\
		\cfl{n}{2} &, \text{ otherwise.}
	\end{cases}$$\\
	This implies that for each value of $p(\offBra{1})$, there exists an integer \(k\) such that $p(\offBra{1})+k=\floor{n}{2}$, where $l\leq k \leq \ccl{n}{2}-1$. In each $\offBra{j}$ $(2\leq j\leq k)$ at $u$, label a pendant vertex(or isolated vertex if the order that branch is 1) and $V(\offBra{1})$ by even integers and label the rest of the vertices in $V(T)$ by odd. Therefore, $\rna{T}\leq \ceil{\Delta}{2}$.
	\\
	\textbf{Case 2.2}: $l=0,$ $4|n$ and $p(\offBra{1})=n/4$.\\
	As $p(\offBra{1})=n/4$, we get $2\leq p(\offBra{2})\leq n/4$. For each value of $p(\offBra{2})$, there exists an integer $0\leq k \leq \frac{n}{2}-2$ such that $p(\offBra{1})+p(\offBra{2})+k=\frac{n}{2}$. In each $\offBra{j}$ $(3\leq j\leq k)$ at $u$, label a pendant vertex(or isolated vertex if the order that branch is 1) and $V(\offBra{1})$ and $ V(\offBra{2})$ by even integers and label the rest of the vertices in $V(T)$ by odd. Hence, $\rna{T}\leq 2+k\leq \ceil{\Delta}{2}$.
\end{proof}

\begin{theorem}
	Let $T$ be a tree of order $n\geq 14$, and $\Delta(T)=\ceil{n}{2}$ with $u$ as the vertex of degree $\Delta$. We shall assume that $p(\offBra{j})\geq p(\offBra{j+1}), 1\leq j<\Delta-1$.\\
	Suppose $3\leq p(\offBra{1})< \cfloor{n}{2}$ and $p(\offBra{j})\leq \cceil{n}{2}-j+1$ for all $ 1<j<\cceil{n}{2}$. Then $\rna{T}\leq \ceil{\Delta}{2}$.
\end{theorem}

\begin{proof}
	Given $n\geq 14$, we have that $\floor{n}{2}-1\geq \cceil{n}{2}$ (from lemma \ref{lemma2.1}). Let $d=\ceil{\Delta}{2}$.
	\begin{claim}\label{claim2.5}
		$\sum_{j=1}^{d}p(\offBra{j})\geq \floor{n}{2}$.
	\end{claim}
	
	\begin{proof}
		As $\Delta=\ceil{n}{2}$, we see that
		\begin{equation*}
			\sum_{j=1}^{\Delta}(p(\offBra{j})-1)=\left\lfloor\frac{n}{2}\right\rfloor-1.
		\end{equation*}
		Also, $3\leq p(\offBra{1})\leq \cfloor{n}{2}-1$, which implies
		\begin{equation*}
			\left\lfloor\frac{n}{2}\right\rfloor-\left\lceil\frac{1}{2}\left\lfloor\frac{n}{2}\right\rfloor\right\rceil+1\leq \sum_{j=2}^{\Delta}(p(\offBra{j})-1)\leq \left\lfloor\frac{n}{2}\right\rfloor-3.
		\end{equation*}
		As $\ceil{n}{2}-1>\floor{n}{2}-3$, there exists an integer $m>0$ such that for all $i>m$, $p(\offBra{i})-1=0$.\\
		If $m\leq d$, then $p(\offBra{j})-1=0$ for all $j>m$, and 
		\begin{equation*}
			\left\lfloor\frac{n}{2}\right\rfloor-\left\lceil\frac{1}{2}\left\lfloor\frac{n}{2}\right\rfloor\right\rceil+1\leq \sum_{j=2}^{d}(p(\offBra{j})-1)\leq \left\lfloor\frac{n}{2}\right\rfloor-3.
		\end{equation*}
		Therefore, 
		\begin{equation*}
			\sum_{j=1}^{d}p(\offBra{j})\geq \left\lfloor\frac{n}{2}\right\rfloor+\left\lceil\frac{1}{2}\left\lceil \frac{n}{2}\right\rceil\right\rceil-1.
		\end{equation*}
		If $m>d$, then $p(\offBra{j})-1\geq 1$ for all $j\leq m$, and 
		\begin{equation*}
			\sum_{j=2}^{d}(p(\offBra{j})-1)\geq d-1.
		\end{equation*}
		Therefore, 
		\begin{equation*}
			\sum_{j=1}^{d}p(\offBra{j})\geq 2\left\lceil\frac{1}{2}\left\lceil \frac{n}{2}\right\rceil\right\rceil+1
		\end{equation*}
	\end{proof}
	From the assumption, we get $p(\offBra{1})+\cceil{n}{2}-1<\floor{n}{2}$. Consider the function $f(x)=\sum_1^{x} p(\offBra{j}) + \cceil{n}{2}-x$, where $x$ is an integer. 
	\begin{equation*}\label{ineq1}
		f(k_2)=\sum_1^{k_2} p(\offBra{j}) + \left\lceil\frac{1}{2}\left\lceil\frac{n}{2}\right\rceil\right\rceil-k_2 \geq \left\lfloor\frac{n}{2}\right\rfloor
	\end{equation*}
	where $k_2=\text{min}\{x\in \mathbb{Z}|f(x)\geq \floor{n}{2}\}\leq d$ (from claim \ref{claim2.5}) and $\sum_1^{k_2 -1} p(\offBra{j}) + \cceil{n}{2}-{k_2+1} < \floor{n}{2}$. Also, 
	\begin{equation*}
		\sum_{i=1}^{k_2}p(\offBra{j})\leq \sum_1^{k_2 -1} p(\offBra{j}) + \left\lceil\frac{1}{2}\left\lceil\frac{n}{2}\right\rceil\right\rceil-{k_2+1} < \left\lfloor\frac{n}{2}\right\rfloor
	\end{equation*}
	There exists another integer $k_3\geq 0$ such that $$\sum_1^{k_2} p(\offBra{j}) + \cceil{n}{2}-k_2-k_3= \floor{n}{2}$$  
	This implies $k_3< \cceil{n}{2}-k_2$, otherwise if $k_3\geq \cceil{n}{2}-k_2$, then $\sum_1^{k_2} p(\offBra{j})\geq \floor{n}{2}$, a contradiction. \\
	Label $V(\offBra{i}),1\leq i\leq k_2$, and a pendant vertex of $T$ in  each $\offBra{k},k_2+1\leq k\leq \cceil{n}{2}-k_3$ by even numbers and label the rest of the vertices in $V(T)$ by odd numbers. Therefore, $\rna{T}\leq \ceil{\Delta}{2}$.
\end{proof}
\begin{figure}[h]
	\centering
	
	\begin{tikzpicture}[scale=0.7,rotate around={180:(0,0)}]
		\draw  (0,0)-- (0,3);
		\draw  (0,3)node[below]{$u_0$}-- (-2,1);
		\draw  (0,3)-- (2,1);
		\draw  (0,3)-- (-3,3);
		\draw  (0,3)-- (3,3);
		\draw  (0,0)-- (-0.87,-1.15);
		\draw  (0,0)-- (-0.37,-1.37);
		\draw  (0,0)-- (0.37,-1.37);
		
		\draw  (0,0)-- (0.87,-1.15);
		\draw  (2,1)-- (1.43,-0.15);
		\draw  (2,1)-- (1.85,-0.19);
		\draw  (2,1)-- (2.25,-0.13);
		\draw  (2,1)-- (2.63,0.09);
		\draw  (3,3)-- (2.45,1.95);
		\draw  (3,3)-- (3,2);
		\draw  (3,3)-- (3.47,2.07);
		\draw  (3,3)-- (3.89,2.33);
		\draw  (-2,1)-- (-2.95,0.55);
		\draw  (-2,1)-- (-2.59,0.17);
		\draw  (-2,1)-- (-2.13,-0.01);
		\draw  (-2,1)-- (-1.67,0.09);
		\draw  (-3,3)-- (-3.95,2.83);
		\draw  (-3,3)-- (-3.89,2.39);
		\draw  (-3,3)-- (-3.45,1.99);
		\draw  (-3,3)-- (-2.77,2.03);
		\begin{scriptsize}
			\draw [fill=black, label=$u_0$] (0,0) circle (2pt);
			\draw [fill=black] (0,3) circle (2pt);
			\draw [fill=black] (-2,1) circle (2pt);
			\draw [fill=black] (2,1) circle (2pt);
			\draw [fill=black] (-3,3) circle (2pt);
			\draw [fill=black] (3,3) circle (2pt);
			\draw [fill=black] (-0.87,-1.15) circle (2pt);
			\draw [fill=black] (-0.37,-1.37) circle (2pt);
			\draw [fill=black] (0.37,-1.37) circle (2pt);
			\draw [fill=black] (0.87,-1.15) circle (2pt);
			\draw [fill=black] (1.43,-0.15) circle (2pt);
			\draw [fill=black] (1.85,-0.19) circle (2pt);
			\draw [fill=black] (2.25,-0.13) circle (2pt);
			\draw [fill=black] (2.63,0.09) circle (2pt);
			\draw [fill=black] (2.45,1.95) circle (2pt);
			\draw [fill=black] (3,2) circle (2pt);
			\draw [fill=black] (3.47,2.07) circle (2pt);
			\draw [fill=black] (3.89,2.33) circle (2pt);
			\draw [fill=black] (-2.95,0.55) circle (2pt);
			\draw [fill=black] (-2.59,0.17) circle (2pt);
			\draw [fill=black] (-2.13,-0.01) circle (2pt);
			\draw [fill=black] (-1.67,0.09) circle (2pt);
			\draw [fill=black] (-3.95,2.83) circle (2pt);
			\draw [fill=black] (-3.89,2.39) circle (2pt);
			\draw [fill=black] (-3.45,1.99) circle (2pt);
			\draw [fill=black] (-2.77,2.03) circle (2pt);
		\end{scriptsize}
	\end{tikzpicture}
	
	\caption{$(5,2)$-tree with $u_0$ as the root vertex}
	\label{fig:lobster}
\end{figure}
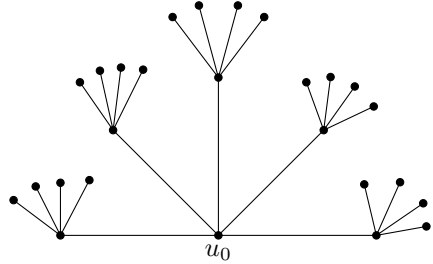
\paragraph{For every tree, $T$, of order $n$, is $\rna{T}\leq \ceil{\Delta(T)}{2}$?}
The answer is no. There exist trees which DO NOT satisfy the above inequality and to prove that, we will define a new class of trees below.
\begin{definition}
	A $(k,l)$-tree is a rooted tree with $u_0$ as the root vertex, in which degree of each non-pendant vertex is $k$ and each pendant vertex is at a distance $l$ from the root vertex $u_0$. Example: figure \ref{fig:lobster}
\end{definition}

\begin{theorem}\label{violeteThm}
	Let $T$ be a $(k,2)-$tree, where $k$ is odd. Then $\rna{T}=k>\ceil{\Delta}{2}$.
\end{theorem}

\begin{proof} 
	Let $k=2m+1,m>0$. The order of $T$ is $k^2+1=4m^2+4m+1$ and $|N_o(k^2+1)|=|N_e(k^2+1)|=2m^2+2m+1$. We can see that degree of each vertex is either 1 or $k$, and the number of pendant vertices is $k(k-1)=4m^2+2m$. If $u_i$ are the neighbours of the root vertex $u_0$, then $deg(u_i)=deg(u_0)=k=2m+1$, for all $1\leq i\leq k$. $N[u_i]-u_0$ induces a star graph, $K_{1,2m}$ (an induced subgraph of $T$), for each $i$, and let us denote these induced subgraphs as $S_i$ and for every $1\leq i\leq k$,  $p(S_i)=k=2m+1$. Also, $\sum_{i=1}^m p(S_i) <2m^2+2m+1$ and $\sum_{i=1}^{m+1} p(S_i) > 2m^2+2m+1$. Let us label $V(S_i)$ for all $i\in N(m)$ and $m+1$ pendant vertices of $T$ in $S_{m+1}$ by even integers and rest of the vertices of $T$ by odd integers. Therefore, $\rna{T}\leq 2m+1$.\\
	We shall prove that $\rna{T}$ is, in fact, equal to $2m+1$.
	Let $A, B$ be a bipartition of $V(T)$ and $E(A,B)$ be the set of edges whose one end vertex is in $A$ and  the other end vertex is in $B$. 
	We shall investigate the influence of $E(A,B)$ (type of edges in $E(A,B)$) on the \textit{rna} number of $T$. For the sake of convenience, we shall define $S_0=K_0$, a null graph or empty graph- a graph with no vertices and edges.\vspace{2mm}
	
	\noindent\textbf{Case 1}: If the edges in $E(A,B)$ are pendant edges of $T$, then $|E(A,B)|=2m^2+2m+1>2m+1$. Therefore, if the edges in $E^-(T,\sigma)$ are pendant edges of $(T,\sigma)$, then $|E^-(T,\sigma)|>2m+1$, under some parity signature $\sigma$.\vspace{2mm}
	
	\noindent\textbf{Case 2}: There are $k$ non-pendant edges in $T$ and each of them is incident with the root vertex $u_0$. Suppose the edges in $E(A,B)$ are non-pendant edges of $T$. It is enough to study when $|E(A,B)|$ is either $\ceil{k}{2}$ or $\floor{k}{2}$ because if $|A|\neq |B|$ whenever $|E(A,B)|$ is either $\ceil{k}{2}$ or $\floor{k}{2}$, then $|A|\neq |B|$ whenever  $|E(A,B)|$ is either $>\ceil{k}{2}$ or $<\floor{k}{2}$. \\
	Suppose $|E(A,B)|=\floor{k}{2}$ with $u_0\in A$. Then $|A|=1+k\ceil{k}{2}=2m^2+3m+2>2m^2+2m+1$ and $|B|=k\floor{k}{2}=2m^2+m<2m^2+2m+1$. Therefore, $|A|\neq |B|$ and $A$ and $B$ does not form a bipartition of $V(T)$.\\
	Suppose $|E(A,B)|=\ceil{k}{2}$ with $u_0\in A$. Then $|A|=1+k\floor{k}{2}=2m^2+m+1<2m^2+2m+1$ and $|B|=k\ceil{k}{2}=2m^2+3m+1>2m^2+2m+1$. Therefore, $|A|\neq |B|$ and $A$ and $B$ does not form a bipartition of $V(T)$.\\
	Therefore, $E(A,B)$ cannot contain only non-pendant edges of $T$.\vspace{2mm}
	
	\noindent\textbf{Case 3}:  Suppose $E(A,B)$ contains both pendant and non-pendant edges. Let $k_1$=number of pendant edges of $T$ in $E(A,B)$ and $k_2=$ number of non-pendant edges of $T$ in $E(A,B)$, and $|E(A,B)|=k_1+k_2$. Also, $k_1\neq 0$ and $k_2\neq 0$. \\
	Without loss of generality, let us assume $u_0\in B$. Let $U_j$ be a subgraph of $T[A]$ induced by $u_j$ and its neighbours in $[A]$ (if $u_j$ has no neighbours in $[A]$, then $U_j$ is an isolated vertex in $T[A]$) such that $1\leq p(U_j)\leq 2m$. Except for $u_0$, all the neighbours of $u_j$ are pendant vertices. If a neighbour, $v\neq u_0$, of $u_j$ is in $T[B]$, then $u_jv\in E(A,B)$. Therefore, $u_j$ is adjacent to $p(S_j)-p(U_j)=2m+1-p(U_j)$ vertices in $T[B]$ and hence, $k_1$ is atleast $\sum_{j}(2m+1-p(U_j))$.\\
	$T[A]$ can have isolated vertices, which are basically the pendant vertices adjacent to some vertices in $T[B]$ and therefore, pendant edges associated with these isolated vertices are present between $T[A]$ and $T[B]$. Either  $\cup_{i=1}^{k}S_i \subseteq T[A]$ or $T[A]$ has no subgraphs $S_i$, where $k\leq m$.\\
	\textbf{Case 3.1}: Suppose $k_2<m$ or $k_2=m-b, b>0$. Then $\{u_i\}_{i=1}^{k_2} \subseteq A$.\\
	\textbf{Case 3.1.1}\label{3.1.1}: Suppose $\cup_{i=1}^{k_2}S_i \subseteq T[A]$.\\ Now, $\sum_{i=1}^{k_2}p(S_i)=k_2(2m+1)=2m^2+m-2mb-b<2m^2+2m+1$ and thus there are $2mb+m+b+1$ isolated vertices in $T[A]$. Therefore, $k_1=2mb+m+b+1$ (no $u_j$ has neighbours in $T[B]$) and 
	\begin{equation*}
		|E(A,B)|=k_1+k_2=2m+2mb+1>2m+1
	\end{equation*}
	\textbf{Case 3.1.2}: Let $0\leq k^\prime<k_2$ be an integer. 
	Now $\sum_{i=0}^{k^\prime}p(S_i)=k^\prime(2m+1)<k_2(2m+1)<2m^2+2m+1$. Further, $1\leq p(U_j)\leq 2m$, for $k^\prime+1\leq j\leq k_2$. Therefore, 
	\begin{equation*}
		\sum_{i=0}^{k^\prime}p(S_i)+\sum_{j=k^\prime+1}^{k_2}p(U_j)<\sum_{l=0}^{k_2}p(S_l)
	\end{equation*}
	Also, there are more than $2mb+m+b+1$ isolated vertices in $T[A]$. Therefore,
	\begin{align*}
		k_1 & >\sum_{j=k^\prime+1}^{k_2}(2m+1-p(U_j))+2mb+m+b+1\\
		|E(A,B)|&=k_1+k_2>2m+1  
	\end{align*} 
	\textbf{Case 3.2}: Suppose $m+1<k_2\leq 2m+1$ or $k_2=m+y$, where $2\leq y\leq m+1$. $T[A]$ can contain at most $m$ induced subgraphs $S_i$. Let $T[A]$ have $z$ induced subgraphs $S_i$, where $0\leq z\leq m$. There are two subcases. \\	
	\textbf{Case 3.2.1}: Suppose $k_2=2m+1$. Then
	\begin{equation*}
		k_1\geq\sum_{j=z+1}^{2m+1}(2m+1-p(U_j))\geq 2m-z+1
	\end{equation*}
	Therefore, $|E(A,B)|=k_1+k_2\geq (2m-z+1)+(2m+1)>2m+1$.\\
	\textbf{Case 3.2.2}: Suppose $k_2<2m+1$ and $\cup_{i=1}^z S_i\subseteq T[A]$. \\
	If $m=1$, then $2m+1=3$ and $m+1=2$, thus, $k_2$ doesn't exist as $2<k_2<3$.\\
	Let $m>1$. Then, $\sum_{j=z+1}^{k_2}p(U_j)\leq (k_2-z)(2m)\leq 4m^2-2mz$. As there are $z$ induced subgraphs $S_i$ in $T[A]$, the sum $\sum_{j=z+1}^{k_2}p(U_j)$ cannot exceed $|A|-\sum_{i=0}^{z}p(S_i)$ and  $$|A|-\sum_{i=0}^{z}p(S_i)\geq (2m^2+m+1)-z(2m+1)=2m^2+2m-2mz-z+1$$
	Furthermore, $2m^2+2m-2mz-z+1<4m^2+2mz$. Therefore, $k_1\geq \sum_{j=z+1}^{k_2}(2m+1-p(U_j))\geq 2m^2-2m-1+z$ and 
	\begin{equation*}
		|E(A,B)|=k_1+k_2> (2m^2-2m-1+z)+(m+1)=2m^2-m+z>2m+1
	\end{equation*}
	\textbf{Case 3.3}: $k_2=m$ and $\cup_{i=1}^{k_3}S_i \subseteq T[A]$, where $0\leq k_3\leq m$. $[A]$ has at least $m+1$ isolated vertices. There are two subcases to consider.\\
	\textbf{Case 3.3.1}: If $k_3=m$, then $T[A]$ has $m+1$ isolated vertices and $k_1=m+1$. Therefore,
	\begin{equation*}
		|E(A,B)|=k_1+k_2=(m+1)+(m)=2m+1
	\end{equation*}
	\textbf{Case 3.3.2}: If $k_3<m$, then $k_1>m+1+\sum_{j=k_3+1}^{m}(2m+1-p(U_j)>m+1$. Therefore, $$|E(A,B)|=k_1+k_2>(m+1)+m=2m+1$$\\
	\textbf{Case 3.4}: $k_2=m+1$ and $\cup_{i=1}^{k_4}S_i \subseteq T[A]$, where $0\leq k_4\leq m$.\\
	\textbf{Case 3.4.1}: If $k_4=m$ and $p(U_{m+1})\leq m+1$, then $k_1\geq m$.
	Therefore, $$|E(A,B)|=k_1+k_2\geq (m)+(m+1)=2m+1$$
	If $p(U_j)=m+1$, then $|E(A,B)|=2m+1$ \\
	\textbf{Case 3.4.2}: If $k_4<m$, then $k_1\geq \sum_{j=k_4+1}^{m+1}(2m+1-p(U_j))\geq (m-k_4+1)(2m)$.
	Therefore, $$|E(A,B)|=k_1+k_2\geq 2m(m-k+1)+(m+1)>3m+1$$ 
	From the above cases we see that, $|E(A,B)|=2m+1$ only in the case 3.3.1 and the case 3.4.1 and in the rest of the cases, $|E(A,B)|>2m+1$. Therefore, $\rna{T}=2m+1>\ceil{\Delta}{2}$, where $\Delta=2m+1$.
\end{proof}
\section{An important property of the \textit{rna} number of trees}\label{sec2.2}

Acharya et al proved that there exists graphs with a desired $rna$ number\cite[Theorem 3.4]{acharya_characterizations_2021}. We ask a similar question in the case of trees: For a given order $n$, are there trees of order $n$ with a desired $rna$ number? This is answered in the next theorem. 

A \textit{spider} is a tree having at most one vertex (called the central vertex) of degree greater than two (\cite{gallian2018dynamic}). We denote a spider graph by $S_m(l_1,l_2,\ldots,l_k)$, where \(m\) is the number of branches (called the legs) at the vertex $u$, which has $k$ branches of length $ l_i>1$, where $1\leq i\leq k$, and $m-k$ branches of length 1. Let us consider spiders with $k=1$, denoted by $S_m(l)$, whose order is $n+1$. For the sake of completeness, we shall denote $S_m(1)=K_{1,n}$.

\begin{theorem}\label{similarthm}
	For a given integer $n+1>1$ and for each $i \in \{1,2,\ldots,\ceil{n}{2}\}$, there exists a tree of order $n+1$ whose \(rna\)-number is $i$.
\end{theorem}
To prove Theorem \ref{similarthm}, it is enough to prove the following statement
\begin{theorem}\label{partition}\label{similarThm}
	$\rna{S_m(l)}=\floor{n+1}{2}-l+1$, $1\leq l \leq \floor{n+1}{2}$, where $n+1>1$.
\end{theorem}

\begin{proof}
	Since $S_m(1)=K_{1,n}$, we know that, $\rna{K_{1,n}}=\ceil{n}{2}=\floor{n+1}{2}$ \cite[Proposition 2.7]{acharya_parity_2021}. We shall consider $S_m(l)$ for $l>1$.
	\begin{claim}
		$\rna{S_m(l)}\leq \begin{cases}
			\frac{n}{2}-l+1 & \text{; if n+1 is odd}\\
			\frac{n+1}{2}-l+1 & \text{; if n+1 is even}
		\end{cases}$.
	\end{claim}
	\begin{proof}
		Let $Br(u,1)=uv u_1 u_2\ldots, u_{l-1}$ be the branch at vertex $u$ of length \(l\) and let $v_1,v_2,\ldots,v_{n-l}$ be the remaining vertices of $S_m(l)$. Consider the path $P_{l}=Br(u,1)-\{u\}=vu_1u_2\ldots u_{l-1}$ and let $A=\{v_1,v_2,\ldots,v_k\}$ and $B=\{v_{k+1},v_{k+2}\ldots v_{n-l}\}$ be sets of pendant vertices, where $k=\floor{n+1}{2}-l$. \\
		Let $f:V(S_m(l))\longrightarrow N(n+1)$ be a bijective function such that $f(V(P_{l})\cup A)=N_e(n+1)$ and $f(\{u\}\cup B)=N_o(n+1)$. Also, $|V(P_l)\cup A|=\floor{n+1}{2}$ and $|\{u\}\cup B|=\ceil{n+1}{2}$. Let $H=(S_m(l),\sigma_f)$ be the parity signed graph under the signature $\sigma_f$. Therefore, $uv, uv_1, uv_2\ldots uv_k$ are the negative edges in $H$, where $k=\floor{n+1}{2}-l$. Hence, 
		\begin{equation}
			\rna{S_m(l)}\leq \begin{cases}
				\frac{n}{2}-l+1 & \text{; if n+1 is odd}\\
				\frac{n+1}{2}-l+1 & \text{; if n+1 is even}
			\end{cases}
		\end{equation} 
		
	\end{proof}
	\begin{claim}\label{claim3.2.2}
		For any parity labelling $f$, $|E^-(S_m(l),\sigma_f)|\geq \floor{n+1}{2}-l+1$.
	\end{claim}
	\begin{proof}
		We shall show that for any parity labelling $f$, $|E^-(S_m(l),\sigma_f)|\geq \floor{n+1}{2}-l+1$. \\
		Let $S_m(l)[\{u,v,v_1,v_2,\ldots,v_{n-l}\}]=K_{1,n-l+1}$ and $S_m(l)[\{u_1,u_2,\ldots,u_{l-1}\}]=P_l$. $vu_1$ is an edge between $K_{1,n-l+1}$ and $P_l$ in $S_m(l)$.\\ 
		Let $G(n+1)$ be the set of all injective functions $g:V(K_{1,n-l+1})\longrightarrow N(n+1)$ and $W(n+1)$ be the set of all injective functions $w:V(P_l)\longrightarrow N(n+1)-g(V(K_{1,n-l+1}))$. We can define a bijective function $h$ using two such injective functions, $g\in G(n+1)$ and $w\in W(n+1)$: Let $h:V(S_m(l))\longrightarrow N(n+1)$ such that $h|_{V(K_{1,n-l+1})}=g$ and $h|_{V(P_l)}=w$ and we can see that $h$ is a bijective function. Similarly, we can obtain two such injective functions from a given bijective function.\\
		\indent For every injective function $g\in G(n+1)$, let $\sigma_g:E(K_{1,n-l+1})\longrightarrow \{-,+\}$ be a function such that, $\sigma_g(yz)=+$ if $g(y)$ and $g(z)$ are of same parity and $\sigma_g(yz)=-$ if $g(y)$ and $g(v)$ are of opposite parity. Similarly,for every injective function $w\in W(n+1)$ , let $\sigma_w:E(P_l)\longrightarrow \{-,+\}$ be a function such that, $\sigma_w(yz)=+$ if $w(y)$ and $w(z)$ are of same parity and $\sigma_w(yz)=-$ if $w(y)$ and $w(v)$ are of opposite parity. Therefore, $\sigma_h:E(S_m(l))\longrightarrow \{-,+\}$ is a function such that $\sigma_h|_{E(K_{1,n-l+1})}=\sigma_g$ and $\sigma_h|_{E(P_l)}=\sigma_w$, and $\sigma_h(vu_1)= + $ or $-$, depending on the value of $g(v)$ and $w(u_1)$.\\ 
		
		\noindent\textbf{Case 1}: $n+1$ is odd. \\
		We know that $1<l\leq n/2$ which implies $n+1>2l-1$, hence, $\frac{n+2}{2}-1<n-l+1$.
		\textbf{Case 1.1}: $g(u)$ is odd, where $g\in G(n+1)$.\\
		Let $|E^-(K_{1,n-l+1},\sigma_g)|$ be the number of negative edges which is equal to the number of even labelled vertices. As $\frac{n+2}{2}-1<n-l+1$, we can use the remaining $\ceil{n+1}{2}-1$ odd integers to label the pendant vertices of $K_{1,n-l+1}$, that is, $g:V(K_{1,n-l+1})\longrightarrow N_o(n+1)\cup N_e(n-2l+2)$ with $g(u)$=odd. Therefore, $|E^-(K_{1,n-l+1})|=\frac{n}{2}-l+1$ and is the smallest number of even labelled vertices in $K_{1,n-l+1}$ under any labelling $g\in G(n+1)$ because $N(n+1)$ has $\ceil{n+1}{2}$ odd integers and to make $|E^-(K_{1,n-l+1})|<n/2-l+1$, we need more than $\ceil{n+1}{2}$ odd labelled vertices in $K_{1,n-l+1}$ under some labelling which is not possible.\\
		After labelling the vertices of $K_{1,n-l+1}$, we are left with only even integers to label the vertices of $P_l$ and  $|E^-(P_l,\sigma_w)|=0$ under the labelling $w\in W(n+1)$. We can assume that $g(v)$=even number and hence, $|E^-(S_m(l),\sigma_h)| = |E^-(K_{1,n-l+1},\sigma_g)|+|E^-(P_l,\sigma_w)|\geq n/2-l+1$ under any bijective function $h$ defined by $g\in G(n+1)$ and $w\in W(n+1)$, where $\sigma_h(vu_1)=-$ and $g(u)$ is odd.\\
		\textbf{Case 1.2}: $g(u)$ is even, where $g\in G(n+1)$.\\
		Let $|E^-(K_{1,n-l+1},\sigma_g)|$ be the number of negative edges which is equal to the number of vertices with odd integers as their labels. As $\frac{n+2}{2}-1<n-l+1$, we can use the remaining $\floor{n+1}{2}-1$ even integers to label the pendant vertices of $K_{1,n-l+1}$, that is, $g:V(K_{1,n-l+1})\longrightarrow N_e(n+1)\cup N_o(n-2l+4)$ with $g(u)$=even. Therefore, $|E^-(K_{1,n-l+1})|=\frac{n}{2}-l+2$ and is the smallest number of odd labelled vertices in $K_{1,n-l+1}$ under any labelling $g\in G(n+1)$ because $N(n+1)$ has $\floor{n+1}{2}$ even integers and to make $|E^-(K_{1,n-l+1})|<n/2-l+2$, we need more than $\floor{n+1}{2}$ even labelled vertices in $K_{1,n-l+1}$ under some labelling which is not possible.\\
		After labelling the vertices of $K_{1,n-l+1}$, we are left with only odd integers to label the vertices of $P_l$ and  $|E^-(P_l,\sigma_w)|=0$ under the labelling $w\in W(n+1)$. We can assume that $g(v)$=odd number and hence, $|E^-(S_m(l),\sigma_h)| = |E^-(K_{1,n-l+1},\sigma_g)|+|E^-(P_l,\sigma_w)|\geq n/2-l+2$ under any bijective function $h$ defined by $g\in G(n+1)$ and $w\in W(n+1)$, where $\sigma_h(vu_1)=-$ and $g(u)$ is even.\\
		Therefore, when $n+1$ is odd and for any parity labelling $f$, we have  $|E^-(S_m(l),\sigma_f)|\geq n/2-l+1$. \\
		
		\noindent\textbf{Case 2}: $n+1$ is even.\\
		We know that $1<l\leq \frac{n+1}{2}$ which implies $n+2>2l-1$, and hence, $\frac{n+1}{2}-1<n-l+1$. Since $|N_o(n+1)|=|N_e(n+1)|$, without loss of generality, we can assume that $g(u)$ is odd, where $g\in G(n+1)$.\\
		Let $|E^-(K_{1,n-l+1},\sigma_g)|$ be the number of negative edges which is equal to the number of even labelled vertices. As $\frac{n+1}{2}-1<n-l+1$, we can use the remaining $\frac{n+1}{2}-1$ odd integers  and appropriate number of even integers to label the pendant vertices of $K_{1,n-l+1}$, that is, $g:V(K_{1,n-l+1})\longrightarrow N_o(n+1)\cup N_e(n-2l+3)$ with $g(u)$=odd integer. Therefore, $|E^-(K_{1,n-l+1})|=\frac{n+1}{2}-l+1$ and is the smallest number of even labelled vertices in $K_{1,n-l+1}$ under any labelling $g\in G(n+1)$ because $N(n+1)$ has $\frac{n+1}{2}$ odd integers and to make $|E^-(K_{1,n-l+1})|<(n+1)/2-l+1$, we need more than $\frac{n+1}{2}$ odd labelled vertices in $K_{1,n-l+1}$ under some labelling which is not possible.\\ After labelling the vertices of $K_{1,n-l+1}$, we are left with only even integers to label the vertices of $P_l$ and  $|E^-(P_l,\sigma_w)|=0$ under the labelling $w\in W(n+1)$. We can assume that $g(v)$=even number and hence, $|E^-(S_m(l),\sigma_h)| = |E^-(K_{1,n-l+1},\sigma_g)|+|E^-(P_l,\sigma_w)|\geq \frac{n+1}{2}-l+1$ under any bijective function $h$ defined by $g\in G(n+1)$ and $w\in W(n+1)$, where $\sigma_h(vu_1)=-$ and $g(u)$ is odd.\\
		Therefore, when $n+1$ is even and for any parity labelling $f$, we have $|E^-(S_m(l),\sigma_f)|\geq \frac{n+1}{2}-l+1$.
	\end{proof}
	From claim \ref{claim3.2.2}, we get $\rna{S_m(l)}\geq \floor{n+1}{2}-l+1$. Therefore, $\rna{S_m(l)}=\floor{n+1}{2}-l+1$ for all $1\leq l\leq \floor{n+1}{2}$.
\end{proof}

\section{Characteristics of graphs associated with the rna number}\label{sec2.3}
Before we move to the characterizations section, we shall prove some results. Let $S(G)$ be the set of spanning trees of a graph \(G\). Let $(G,\sigma_{f_{rna}})$ be the parity signed graph with a signature $\sigma_{f_{rna}}$ (obtained by an appropriate labelling function $f_{rna}$), which gives the minimum number of negative edges.

\begin{theorem}\label{rna1thm}
	If $\rna{G}=1$, then $\rna{T}=1$ for all $T\in S(G)$.
\end{theorem}

\begin{figure}
	\centering
	\subfigure[]{
		\begin{tikzpicture}[]
			\draw (-2,0)--(-1,0)--(0,1)--(1,0)--(2,0);
			\draw (-1,0)--(1,0);
			\draw[fill=black] (0,1) circle (2pt);
			\draw[fill=black] (-1,0) circle (2pt);
			\draw[fill=black] (-2,0) circle (2pt);
			\draw[fill=black] (1,0) circle (2pt);
			\draw[fill=black] (2,0) circle (2pt);
		\end{tikzpicture} \label{fig2.a}
	}  
	\subfigure[]{
		\begin{tikzpicture}[]
			\draw[fill=black] (0,1) circle (2pt);
			\draw[fill=black] (-1,0) circle (2pt);
			\draw[fill=black] (-2,0) circle (2pt);
			\draw[fill=black] (1,0) circle (2pt);
			\draw[fill=black] (2,0) circle (2pt);
			\draw (-2,0)--(-1,0)--(0,1)--(1,0)--(2,0);
		\end{tikzpicture}
	}
	\subfigure[]{
		\begin{tikzpicture}[]
			\draw[fill=black] (0,1) circle (2pt);
			\draw[fill=black] (-1,0) circle (2pt);
			\draw[fill=black] (-2,0) circle (2pt);
			\draw[fill=black] (1,0) circle (2pt);
			\draw[fill=black] (2,0) circle (2pt);
			\draw (-2,0)--(-1,0)--(1,0)--(2,0);
			\draw (-1,0)--(0,1);
		\end{tikzpicture}
	}
	\caption{(a): A graph with rna number=2. (b) and (c): Spanning trees of the graph (a) with the rna number=1}
\end{figure}
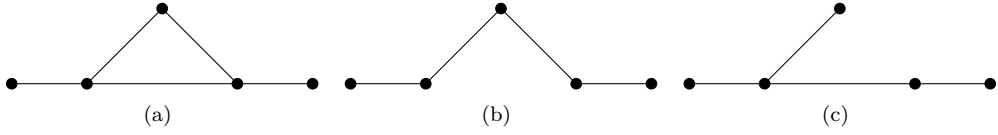

\begin{proof}
	From the assumption, $E^-(G,\sigma_{f_{rna}})=\{e\}$. If $T$ is a spanning tree of \(G\), then $e\in E(T)$ and hence, $E^-(T,\sigma_{f_{rna}})=\{e\}$ under the same labelling ($f_{rna}$) and signature ($\sigma_{f_{rna}}$).  
\end{proof}
The converse is not true. Consider the graph in figure \ref{fig2.a} and its spanning trees in figures 2(b) and 2(c).\\
We will now give a sufficient condition for trees to achieve \(rna\) number=2.
\begin{theorem}
	Let \(T\) be a tree of order $n$. Suppose $C(T)>1$. Let $uv$ be an edge such that $G_1, G_2$ are two components of $T-uv$ and $p(G_2)-p(G_1)=C(T)$. If $T$ satisfies either 
	\begin{enumerate}
		\item $C\left (G_{2},\left\lceil\frac{C(T)}{2}\right\rceil\right)=\left\lfloor\frac{n}{2}\right\rfloor- \left\lceil\frac{C(T)}{2}\right\rceil$ or 
		\item $C \left(G_{2},\left\lfloor\frac{C(T)}{2}\right\rfloor\right)=\left\lceil\frac{n}{2}\right\rceil-\left\lfloor\frac{C(T)}{2}\right\rfloor$
	\end{enumerate}
	then $\rna{T}=2$.
\end{theorem}

\begin{proof}
	Given $G_1, G_2$ are the two components of $T-uv$, and $p(G_2)>p(G_1)$ and 
	\begin{equation}\label{eq0}
		p(G_2)-p(G_1)=C(T)
	\end{equation}
	In this proof, we will prove the statement when $n$ is odd as the same arguments would follow when $n$ is even. As we assumed $n$ to be odd, we get that $C(T)$ is also odd.\\
	(1): Suppose 
	\begin{equation}\label{eq1}
		C\left (G_{2},\left\lceil\frac{C(T)}{2}\right\rceil\right)=\left\lfloor\frac{n}{2}\right\rfloor- \left\lceil\frac{C(T)}{2}\right\rceil
	\end{equation}
	Thus there exists an edge $e_1$ in $G_2$ such that $p(J_1)-p(J_2)=\floor{n}{2} - \ceil{C(T)}{2}$ and $p(J_2)\leq  \ceil{C(T)}{2}$, where $J_1,J_2$ are the two components of $G_2-e$ and $p(J_1)>p(J_2)$. If $p(J_2)=\ceil{C(T)}{2}-b, 0\leq b$, then from equation \eqref{eq1}, $p(J_1)=\floor{n}{2}-b$. $p(G_1)+p(G_2)=n=p(G_1)+\floor{n}{2}+\ceil{C(T)}{2}-2b$ which implies
	\begin{equation}\label{eq1.1}
		p(G_1)=\left\lceil\frac{n}{2}\right\rceil-\left\lceil\frac{C(T)}{2}\right\rceil+2b
	\end{equation}
	From equations \eqref{eq1.1} and \eqref{eq0}, $b=0$. Therefore, we see that $p(G_1)+p(J_2)=\ceil{n}{2}$ and $p(J_1)=\floor{n}{2}$ and hence, $\rna{T}=2$.
	hence $\rna{T}=2$.\\
	(\(2\)): Suppose 
	\begin{equation}\label{eq2}
		C \left(G_{2},\left\lfloor\frac{C(T)}{2}\right\rfloor\right)=\left\lceil\frac{n}{2}\right\rceil-\left\lfloor\frac{C(T)}{2}\right\rfloor
	\end{equation}
	Thus there exists an edge $e_2$ in $G_2$ such that $p(F_1)-p(F_2)=\ceil{n}{2}-\floor{C(T)}{2}$ and $p(F_2)\leq \floor{C(T)}{2} $, where $F_1,F_2$ are the two components of $G_2-f$ and $p(F_1)>p(F_2)$. If $p(F_2)=\floor{C(T)}{2}-b^\prime, 0\leq b^\prime$, then from equation \eqref{eq2}, $p(F_1)=\ceil{n}{2}-b^\prime$. $p(G_1)+p(G_2)=n=p(G_1)+\ceil{n}{2}+\floor{C(T)}{2}-2b^\prime$ which implies
	\begin{equation}\label{eq2.1}
		p(G_1)=\left\lfloor\frac{n}{2}\right\rfloor-\left\lfloor\frac{C(T)}{2}\right\rfloor+2b^\prime
	\end{equation}
	From the equation \eqref{eq2.1} and \eqref{eq0}, $b^\prime=0$. Therefore, we see that, $p(G_1)+p(F_2)=\floor{n}{2}$ and $p(F_1)=\ceil{n}{2}$ and hence, $\rna{T}=2$.
\end{proof}
\subsection{Characterizations}

We will now provide a characterization of graphs with the rna number equal to 1 in terms of the spanning trees.\\

\begin{theorem}
	Let \(G\) be a simple graph and $S(G)$ be the set of all spanning trees, $T$, of \(G\). Let $f_{rna}$ and $\sigma_{f_{rna}}$ be the labelling function and signature, respectively, such that $E^-(G,\sigma_{f_{rna}})=\rna{G}$. Then $\rna{G}=1\text{ if and only if }$ there exists an edge $e \in \bigcap_{T\in S(G)}E(T)$ such that $E^-(T,\sigma_{rna})=\{e\}$ for every $T\in S(G)$  .
\end{theorem}

\begin{proof}
	If $\rna{G}=1$, then there exists an edge $e\in E(G)$, a bijective function $f_{rna}:V(G)\rightarrow N$ and a corresponding signature $\sigma_{f_{rna}}$ such that $E^-(G,\sigma_{f_{rna}})=\{e\}$. From theorem \ref{rna1thm}, $\rna{T}=1$ and for every $T\in S(G)$, sub-signed graph $(T,\sigma_{f_{rna}})$ has only one negative edge, $e$. Therefore, $e\in \bigcap_{T\in S(G)}E(T)$.\\
	Let us prove the sufficient condition. We know that $\rna{T}=1$ for all $T\in S(G)$ and all the signed spanning trees of $G$ share a common negative edge, $e$, under the signature $\sigma_{f_{rna}}$. Therefore, the edge $e$ is a cut-edge in $G$ and $\rna{G}=1$. 
\end{proof}

\indent Two edges $e_1$ and $e_2$ are incident edges if both the edges share the same vertex; otherwise, $e_1$ and $e_2$ are said to be non-incident edges. If the \(rna\)- number of a graph is two, then the two negative edges, $e_1,e_2$, can be of the following three types:
\begin{enumerate}
	\item $e_1,e_2$ belong to a cycle.
	\item $e_1,e_2$ are two incident bridges.
	\item $e_1, e_2$ are two non-incident bridges.
\end{enumerate}
Let $\{G_i\}_i$ be the set of all connected components of $G-E(H)$, where $H$ is a subgraph of \(G\). 
\begin{theorem}
	Let \(G\) be a simple graph of order \(n\) and $C_m=(u_{i})_{i=0}^{m-1}$ be an induced  cycle of length $m$ such that $E(C_m)$ is an edge cut of G. Let $V_1,V_2$ be a parity partition of $V(G)$ such that $|E(V_1,V_2)|=\rna{G}$. Let $\{G_i\}_{i=1}^{k}$ be the set of all connected components of $G-E(C_m)$. The following statements are equivalent.
	\begin{enumerate}
		\item $\rna{G}=2$ with $E(V_1,V_2)\subset E(C_m)$.
		\item \begin{enumerate}
			\item $p(G_i)\leq \ceil{n}{2}$ for all $i\leq k$.
			\item If, for some $j \geq0$, the edge $u_j v_i \in E(G_i)$ is a cut-edge of $G$, then either $deg(u_j)>3$ or $p(G_i)\leq \floor{n}{2}$, where $v_i\in V(G_i)$.
			\item There exists $B\subseteq \{G_i\}_{i=1}^{k}$ such that $H=[V(B)]$ is a connected subgraph of \(G\) of order $\ceil{n}{2}$ or $\floor{n}{2}$, and there exists an integer $j$, $0\leq j\leq m-1$  such that $\{u_{j+i}\}_{i=0}^a\subseteq V(H) \text{ and } V(C_m)\backslash \{u_{j+i}\}_{i=0}^a\subseteq V(G)\backslash V(H)$, where $a\geq 0, j+i\in \mathbb{Z}_m$ which is a set of integers under addition modulo $m$. 
		\end{enumerate}       
	\end{enumerate}
\end{theorem}

\begin{proof}
	$(1)\rightarrow (2)$:\\
	Let $E(V_1,V_2)=\{e_1,e_2\}$, $e_1, e_2\in E(C_m)$ and $G-\{e_1,e_2\}$ has two components $H,K$ such that $p(H),p(K)\in \{\floor{n}{2},\ceil{n}{2}\}$. Therefore, $p(G_i)\leq \ceil{n}{2}$ for all $i\leq k$ and each $G_i$ is a subgraph of either $H$ or $K$. This implies that there exists $B\subseteq \{G_i\}_{i=1}^{k}$ such that $H=[V(B)]$. \\
	We know that $C_m-\{e_1,e_2\}$ is a disjoint union of two paths $P_s, P_t$ (let us consider isolated vertex as a path with one vertex $P_1$) and we can assume that $P_s$ is a subgraph of $H$ and $P_t$ is a subgraph of $K$. Therefore, there exists an integer $a\geq 0$ and an integer $j$ ($0\leq j\leq m-1$) such that $P_s=u_{j}u_{j+1},\ldots, u_{j+a}$ and $j+i\in \mathbb{Z}_m, 0\leq i\leq q$. Furthermore, suppose that an edge $u_j v_i\in E(G_i)$ is a cut-edge of $G$. If $deg(u_j)=3$, then $p(G_i)\leq \floor{n}{2}$, otherwise $\rna{G}=1$ with $E(V_1,V_2)=\{u_j v\}$ (a contradiction) and similarly, if $p(G_i)>\floor{n}{2}$, then $deg(u_j)>3$.\\
	\textbf{$(2)\rightarrow (1)$}:\\
	By conditions $2.(a)\text{ and }2.(c)$, $\rna{G}\leq 2$, and by the condition $2.(b)$, $\rna{G}\neq 1$.
\end{proof}

\begin{theorem}
	Let $uv,vw$ be two incident cut-edges in a graph \(G\) and $deg(v)\geq 3$. Let $G_1,G_2$ and $G_3$ be three components of $G-\{uv,vw\}$, where $v\in V(G_3)$. If $V_1,V_2$ is a parity partition of $G$ such that $|E(V_1,V_2)|=\rna{G}$, then the following are equivalent.
	\begin{enumerate}
		\item $\rna{G}=2$ with $E(V_1,V_2)= \{uv,vw\}$
		\item \begin{enumerate}
			\item $\floor{n}{2}\leq p(G_1)+p(G_2)\leq \ceil{n}{2}$.
			\item $p(G_1),p(G_2)<\floor{n}{2}$.
			\item If $p(G)$ is odd and $p(G_3)=\ceil{n}{2}$, then $deg(v)>3$. 
		\end{enumerate}
	\end{enumerate}
\end{theorem}
\begin{proof}
	$(1)\rightarrow (2)$:\\
	$\rna{G}=2$ and $E(V_1,V_2)=\{uv,vw\}$ implies that $p(G_1)+p(G_2)=\ceil{n}{2} \text{or} \floor{n}{2}$ and that $p(G_1),p(G_2)<\floor{n}{2}$. Suppose $p(G)$ is odd and $p(G_3)=\ceil{n}{2}$. If $deg(v)=3$, then the edge $vx$, $x\in N(v)-\{u,w\}$, is a cut-edge of G and $\rna{G}=1$.\\
	$(2)\rightarrow(1)$:\\
	When $p(G_3)=\floor{n}{2}$, due to the conditions $2.(a) \text{ and } 2.(b)$, $E(V_1,V_2)=\{uv,vw\}$. And from the condition $2.(c)$, if, for any $x\in N(v)-\{u,w\}$, $vx$ is a cut-edge of $G$, then $E(V_1,V_2)\neq\{vx\}$.
\end{proof}
\begin{theorem}
	Let \(G\) be a graph with two non-incident bridges $uv$ and $xy$, and $G_1,G_2$ and $G_3$ be three components of $G-\{uv,xy\}$, where $\{v,y\}\in V(G_3)$. If $V_1,V_2$ is a parity partition of $G$ such that $E(V_1,V_2)=\rna{G}$, then the following are equivalent.
	\begin{enumerate}
		\item $\rna{G}=2$ with $E(V_1,V_2)=\{uv,xy\}$.
		\item \begin{enumerate}
			\item $\floor{n}{2}\leq p(G_1)+p(G_2)\leq \ceil{n}{2}$.
			\item $p(G_1),p(G_2)<\floor{n}{2}$.
			\item Suppose $e$ is a bridge in $G_3$. If $H,K$ are two components of $G-e$, then either $p(H)<\floor{n}{2}$ or $p(K)<\floor{n}{2}$.
		\end{enumerate}
	\end{enumerate}
\end{theorem}

\begin{proof}
	$(1)\rightarrow (2)$:\\
	It is easy to see that $p(G_1)+p(G_2)=\ceil{n}{2} \text{or} \floor{n}{2}$, and that both  $p(G_1),p(G_2)<\floor{n}{2}$. Suppose there is a bridge $e$ in $G_3$ with $H$ and $K$ as the components of $G-e$. Then either $p(H)<\floor{n}{2}$ or $p(K)<\floor{n}{2}$ as $\rna{G}=2$.\\
	$(2)\rightarrow(1)$:\\
	Conditions $2.(a)$ and $2.(b)$ imply that $\rna{G}\leq 2$ and condition $2.(c)$ imply that $\rna{G}\neq1$.
	
\end{proof}
 \bibliographystyle{ieeetr} 
 \bibliography{main.bib}

\end{document}